\documentclass{amsart}

\usepackage{amssymb}
\usepackage{graphicx}
\usepackage{enumerate}
\usepackage{amscd}

\title{A note on sections of broken Lefschetz fibrations}
\author{Kenta Hayano}
\address{Department of Mathematics, Graduate School of Science, Osaka University, Toyonaka, Osaka 560-0043, Japan}
\email{k-hayano@cr.math.sci.osaka-u.ac.jp}

\theoremstyle{plain}
\newtheorem{thm}{Theorem}[section]
\newtheorem{cor}[thm]{Corollary}
\newtheorem{lem}[thm]{Lemma}
\newtheorem{prop}[thm]{Proposition}

\theoremstyle{definition}
\newtheorem{defn}[thm]{Definition}
\newtheorem{exmp}[thm]{Example}

\newtheorem{rem}[thm]{Remark}

\makeatletter

\@addtoreset{figure}{section}
\makeatother

\makeatletter

\@addtoreset{equation}{section}
\makeatother

\begin{document}

\maketitle

\begin{abstract}

We show that there exists a non-trivial simplified broken Lefschetz fibration which has infinitely many homotopy classes of sections. 
We also construct a non-trivial simplified broken Lefschetz fibration which has a section with non-negative square. 
It is known that no Lefschetz fibration satisfies either of the above conditions. 
Smith proved that every Lefschetz fibration has only finitely many homotopy classes of sections, 
and Smith and Stipsicz independently proved that a Lefschetz fibration is trivial if it has a section with non-negative square. 
So our results indicate that there are no generalizations of the above results to broken Lefschetz fibrations. 
We also give a necessary and sufficient condition for the total space of a simplified broken Lefschetz fibration with a section admitting a spin structure, 
which is a generalization of Stipsicz's result on Lefschetz fibrations. 

\end{abstract}

\section{Introduction}

A broken Lefschetz fibration is a smooth map from a $4$-manifold to a $2$-manifold 
which has at most two types of singularities, that is, Lefschetz singularity and indefinite fold singularity. 
Such a fibration was first introduced in \cite{ADK} as a generalization of Lefschetz fibrations to near-symplectic setting. 
Broken Lefschetz fibrations have properties similar to those of Lefschetz fibrations in some aspects. 
So it is natural to try to study the former by the techniques used to develop the latter 
and some of such attempts were successful (e.g. \cite{BK}, \cite{H} and \cite{Per}). 

On the other hand, there are also some crucial differences between two kinds of fibrations. 
For example, it is proved in \cite{AK}, \cite{Ba2} and \cite{Lek} that every closed oriented smooth $4$-manifold admits a broken Lefschetz fibration 
(furthermore, we can prove by using the results in \cite{Lek} and \cite{Wil} that 
every closed oriented smooth $4$-manifold admits a {\it simplified} broken Lefschetz fibration). 
However, there exist a lot of $4$-manifolds which never admits any Lefschetz fibrations 
since the total space of a Lefschetz fibration is symplectic \cite{Gom}. 
So it is important to study how far broken Lefschetz fibrations are different from Lefschetz fibrations. 

Smith proved the following theorem as a generalization of Manin's theorem. 

\begin{thm}[Smith \cite{Smi}]

Let $f:M\rightarrow S^2$ be a non-trivial relatively minimal Lefschetz fibration. 
Then $f$ has only finitely many homotopy classes of sections. 

\end{thm}

The following result implies that we cannot generalize Smith's result to simplified broken Lefschetz fibrations. 

\begin{thm}\label{main1}

For any $g\geq 2$, there exists a non-trivial genus-$g$ simplified broken Lefschetz fibration $f:M\rightarrow S^2$ such that 
no fiber of $f$ contains $(-1)$-sphere and $f$ has infinitely many homotopy classes of sections. 

\end{thm}

Smith and Stipsicz found a constraint on self-intersection numbers of sections of Lefschetz fibrations. 

\begin{thm}[Smith \cite{Smi}, Stipsicz \cite{Sti}]

Let $f:M\rightarrow S^2$ be a genus-$g$ relatively minimal Lefschetz fibration ($g\geq 2$). 
If $f$ has a section $\sigma:S^2\rightarrow M$ which satisfies $[\sigma(S^2)]^2\geq 0$, 
then $f$ is trivial. 

\end{thm}

The following result indicates existence of non-trivial simplified broken Lefschetz fibrations having a section with non-negative square. 

\begin{thm}\label{main2}

For any integer $n\in\mathbb{Z}$ and $g\geq 2$, there exists a non-trivial genus-$g$ simplified broken Lefschetz fibration $f:M\rightarrow S^2$ such that 
$f$ has a section $\sigma:S^2\rightarrow M$ with $[\sigma(S^2)]^2=n$. 

\end{thm}

\begin{rem}

Baykur had already proved in \cite{Ba} that there exists a non-trivial {\it genus-$1$} simplified broken Lefschetz fibration 
which has a section with positive square. 

\end{rem}

In section \ref{pre}, we give a precise definition of (simplified) broken Lefschetz fibrations and 
review monodromy representations of simplified broken Lefschetz fibrations. 
This representation relates the structure of the fibrations to mapping class groups of closed surfaces. 

In section \ref{pro1}, we prove Theorem \ref{main1} by using Kirby diagrams. 
We construct a simplified broken Lefschetz fibration and a family of sections of it. 
To prove any two sections in the family are not homotopic, we calculate the second homotopy group of the total space of the fibration. 

In section \ref{pro2}, we prove Theorem \ref{main2} 
after proving a certain lemma about the relation between monodromy representations and self-intersection numbers of sections. 

In section \ref{spin}, we give a necessary and sufficient condition for the total space of a simplified broken Lefschetz fibration 
admitting a spin structure (Theorem \ref{main3}). 
This result is a generalization of Stipsicz's result \cite{Sti2} on spin structures of total spaces of Lefschetz fibrations. 
After making some remarks about Theorem \ref{main3}, we give some applications of the theorem. 
We prove that the total spaces of some of the fibrations constructed in the proof of Theorem \ref{main2} 
admit spin structures. 
We also give a partial answer to Conjecture 5.3 in \cite{H}, 
which is a conjecture on classification of total spaces of genus-$1$ simplified broken Lefschetz fibrations. 
We prove that the conjecture is true under the hypothesis that the total space of the fibration is spin.

\section{Preliminaries}\label{pre}

\subsection{Broken Lefschetz fibrations}

\begin{defn}

Let $M$ and $B$ be compact oriented smooth manifolds of dimension $4$ and $2$, respectively. 
A smooth map $f:M\rightarrow B$ is called a {\it broken Lefschetz fibration} if it satisfies the following conditions: 

\begin{enumerate}[(1)]

\item $\partial M=f^{-1}(\partial B)$; 

\item $f$ has at most the following types of singularities: 

\begin{itemize}

\item $(z_1,z_2)\mapsto \xi=z_1z_2$, 
where $(z_1,z_2)$ (resp. $\xi$) is a complex local coordinate of $M$ (resp. $B$) compatible with its orientation; 

\item $(t,x_1,x_2,x_3)\mapsto (y_1,y_2)=(t,{x_1}^2+{x_2}^2-{x_3}^2)$, 
where $(t,x_1,x_2,x_3)$ (resp. $(y_1,y_2)$) is a real local coordinate of $M$ (resp. $B$). 

\end{itemize}

\end{enumerate}

The singularities in the condition (2) of the definition are called a {\it Lefschetz singularity} and an {\it indefinite fold singularity}, respectively. 
For a broken Lefschetz fibration $f$, we denote by $\mathcal{C}_f$ (resp. $Z_f$) the set of Lefschetz singularities (resp. indefinite fold singularities) of $f$. 
We call $f$ a {\it Lefschetz fibration} if $Z_f=\emptyset$. 

\end{defn}

In this paper, we will call broken Lefschetz fibrations (resp. Lefschetz fibrations) BLF (resp. LF), for short. 

Let $f:M\rightarrow S^2$ be a BLF. 
We assume that the restriction of $f$ to the set of singularities is injective, $Z_f$ is connected and all the fibers of $f$ are connected. 
Then the set $Z_f$ is either the empty set or an embedded circle in $M$. 
If $Z_f$ is empty, $f$ is an LF over $S^2$. 
If $Z_f$ is an embedded circle, the image $f(Z_f)$ divides the target $2$-sphere into two $2$-disks. 
We denote by $\nu f(Z_f)$ a tubular neighborhood of $f(Z_f)$ and we put
\[
S^2\setminus\text{int}\nu f(Z_f)=D_1\coprod D_2, 
\]
where $D_1$ and $D_2$ are $2$-disks. 
It is easy to see that the genus of a regular fiber of the fibration $f:f^{-1}(D_i)\rightarrow D_i$ is just one higher than that of $f:f^{-1}(D_j)\rightarrow D_j$. 
We call $f^{-1}(D_i)$ (resp. $f^{-1}(D_j)$) the {\it higher side} (resp. {\it lower side}) of $f$ and $f^{-1}(\nu f(Z_f))$ the {\it round cobordism} of $f$. 

\begin{defn}

A BLF $f:M\rightarrow S^2$ is said to be {\it simplified} if it satisfies the following conditions: 

\begin{enumerate}[(1)]

\item $f|_{Z_f\cup \mathcal{C}_f}$ is injective; 

\item $Z_f$ is connected and all the fibers of $f$ are connected; 

\item If $Z_f$ is not empty, $\mathcal{C}_f$ is contained in the higher side of $f$. 

\end{enumerate}

\noindent
For a simplified BLF $f$, the genus of a regular fiber in the higher side of $f$ is called the {\it genus} of $f$. 

\end{defn}

The following lemma was proved by Baykur \cite{Ba}: 

\begin{lem}[Baykur \cite{Ba}]\label{roundhandle}

Let $f$ be a simplified BLF and we denote the higher side and the round cobordism of $f$ by $M_h$ and $M_r$, respectively. 
Then $M_h\cup M_r$ is obtained by $2$-handle attachment to $M_h$ followed by $3$-handle attachment. 
Moreover, the attaching circle of the $2$-handle is a non-separating simple closed curve in a regular fiber of $\text{res}f:M_h\rightarrow D^2$ 
and the framing of the $2$-handle is along the regular fiber. 

\end{lem}

We call an attaching circle of the $2$-handle in the above lemma a {\it vanishing cycle} of the indefinite fold of $f$.

\subsection{Monodromy representations}

Let $f:M\rightarrow B$ be a genus-$g$ LF and $\mathcal{C}_f$ the set of Lefschetz singularities. 
We fix a point $y_0\in B\setminus f(\mathcal{C}_f)$. 
Then a certain homomorphism $\varrho_f:\pi_1(B\setminus f(\mathcal{C}_f),y_0)\rightarrow \mathcal{M}_g$, called a {\it monodromy representation} of $f$, is defined, 
where $\mathcal{M}_g=\text{Diff}^+\Sigma_g/\text{Diff}^+_0\Sigma_g$ is the mapping class group of the genus-$g$ closed oriented surface. 
(for the precise definition of this homomorphism, see \cite{GS}). 

We assume that $B=D^2$ and we put $f(\mathcal{C}_f)=\{y_1,\ldots,y_n\}$. 
We take embedded paths $\alpha_1,\ldots,\alpha_n$ in $D^2$ satisfying the following conditions: 

\begin{itemize}

\item each $\alpha_i$ connects $y_0$ to $y_i$; 

\item if $i\neq j$, then $\alpha_i\cap\alpha_j=\{y_0\}$; 

\item $\alpha_1,\ldots,\alpha_n$ appear in this order when we travel counterclockwise around $y_0$. 

\end{itemize}

We obtain $a_i\in \pi_1(D^2\setminus\{y_1,\ldots,y_n\},y_0)$ ($i=1,\ldots,n$) by connecting a counterclockwise circle around $y_i$ to $y_0$ by using $\alpha_i$. 
We put $W_f=(\varrho_f(a_1),\ldots,\varrho_f(a_n))\in{\mathcal{M}_g}^n$. 
This sequence is called a {\it Hurwitz system} of $f$. 
Kas proved in \cite{Kas} that each $\varrho_f(a_i)$ is the right-handed Dehn twist along a simple closed curve $c_i$ in $\Sigma_g$. 
$c_i$ is called a {\it vanishing cycle} of $y_i$. 

Let $f:M\rightarrow S^2$ be a simplified BLF with $Z_f\neq \emptyset$ and $M_h$ the higher side of $f$. 
Then the restriction of $f$ to $M_h$ is an LF over $D^2$. 
So the monodromy representation and a Hurwitz system of this LF can be defined 
and are called the {\it monodromy representation} and a {\it Hurwitz system} of $f$, respectively. 

\begin{lem}[Baykur \cite{Ba}]

Let $f:M\rightarrow S^2$ be a simplified BLF and $\varrho_f$ a monodromy representation of $f$. 
Then a vanishing cycle $c$ of the indefinite fold of $f$ is preserved by $\varrho_f([\partial D^2])$ up to isotopy. 

\end{lem}

We denote by $\mathcal{M}_g(\gamma)$ the subgroup of $\mathcal{M}_g$ which consists of elements represented by maps preserving the simple closed curve $\gamma$ in $\Sigma_g$ up to isotopy. 
The above lemma says that $\varrho_f([\partial D^2])$ is in $\mathcal{M}_g(c)$ for a vanishing cycle $c$ of the indefinite fold of $f$. 
There is a natural homomorphism $\varphi_c:\mathcal{M}_g(c)\rightarrow \mathcal{M}_{g-1}$ defined by cutting the surface $\Sigma_g$ along $c$ and pasting two $2$-disks along the boundary. 

\begin{lem}[Baykur \cite{Ba}]

The element $\varrho_f([\partial D^2])$ is in the kernel of $\varphi_c$. 
Conversely, for a sequence of simple closed curves $c,c_1,\ldots,c_n$ in $\Sigma_g$ satisfying $t_{c_1}\cdot\cdots\cdot t_{c_n}\in\text{Ker}\varphi_c$, 
there exists a simplified BLF $f:M\rightarrow S^2$ such that a Hurwitz system of $f$ is $(t_{c_1},\ldots,t_{c_n})$ 
and a vanishing cycle of the indefinite fold of $f$ is $c$. 

\end{lem}

\begin{rem}

Such a simplified BLF $f$ is not unique even up to diffeomorphism of the total space. 
Indeed, there exist infinitely many simplified BLFs such that Hurwitz systems of these fibrations are all equivalent 
but the total spaces of these fibrations are mutually not diffeomorphic (see \cite{BK} or \cite{H}). 

\end{rem}

\section{Infinitely many homotopy classes of sections}\label{pro1}

To prove Theorem \ref{main1}, we first give genus-$g$ simplified BLF $f_g:M_g\rightarrow S^2$ and look at the set $[S^2,M_g]$. 
We then construct a family of its sections and prove that any two sections in the family are not homotopic. 

\begin{figure}[htbp]
\begin{center}
\includegraphics[width=145mm]{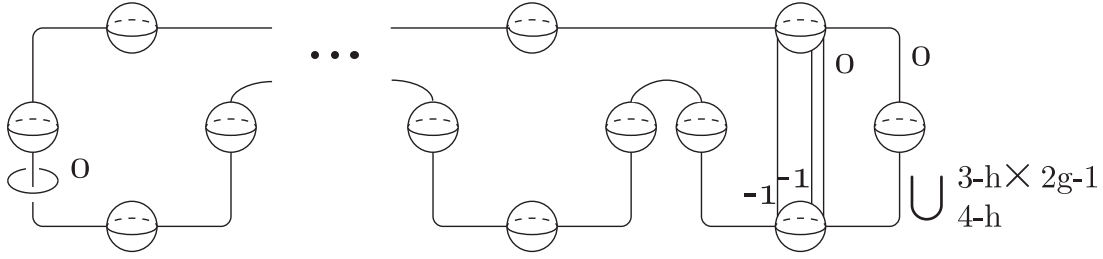}
\end{center}
\caption{The diagram of the total space $M_g$ of $f_g$. 
$2g$ $1$-handles are included in this diagram. }
\label{homotopySBLF}
\end{figure}

\begin{figure}[htbp]
\begin{center}
\includegraphics[width=45mm]{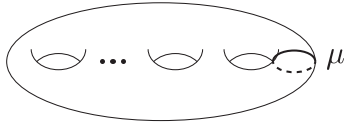}
\caption{The simple closed curve $\mu$ in $\Sigma_g$. }
\label{vanishing}
\end{center}
\end{figure}

\noindent
{\it (Proof of Theorem \ref{main1})}: 
For $g\geq 2$, we denote by $f_g:M_g\rightarrow S^2$ a simplified BLF as shown in Figure \ref{homotopySBLF}. 
This diagram describes the total space of a simplified BLF whose Hurwitz system is $(t_\mu,t_\mu)$, 
where $\mu\subset\Sigma_g$ is a simple closed curve described in Figure \ref{vanishing}. 

\begin{figure}[htbp]
\begin{center}
\includegraphics[width=145mm]{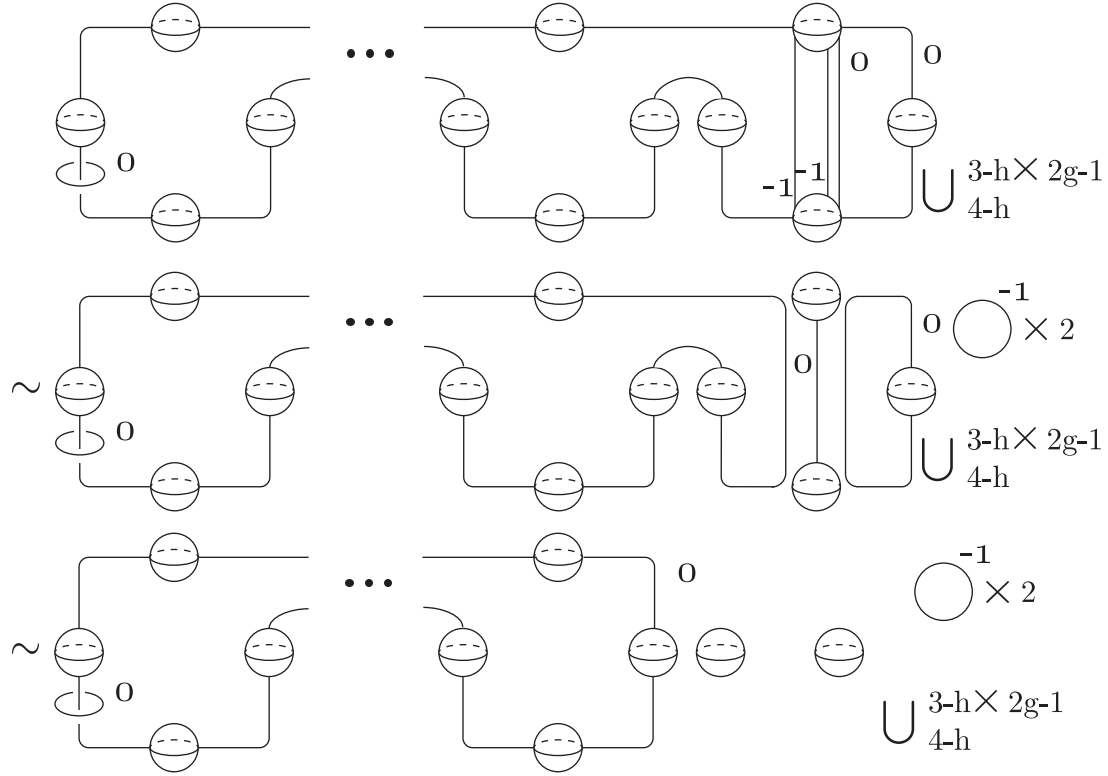}
\caption{The diagram of $M_g$. }
\label{moveSBLF}
\end{center}
\end{figure}

We can change the diagram of $M_g$ as shown in Figure \ref{moveSBLF}, and obtain: 
\[
M_g\cong S^2\times\Sigma_{g-1}\sharp S^1\times S^3 \sharp 2\overline{\mathbb{CP}^2}. 
\]

\begin{figure}[htbp]
\begin{center}
\includegraphics[width=145mm]{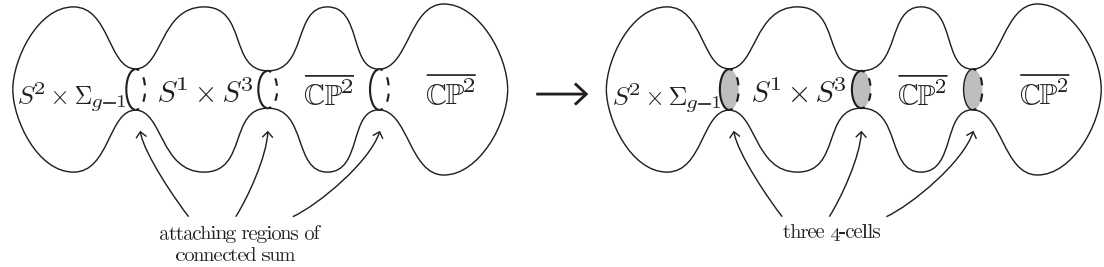}
\caption{Left: the figure describing $M_g$. 
Right: the figure describing $X_g$. 
The shaded parts represent the attached $4$-cells. }
\label{CWcpx}
\end{center}
\end{figure}

To analyze the set $[S^2,M_g]$, we first look at the group $\pi_2(M_g,p_0)$ for a fixed point $p_0\in M_g$. 
Let $X_g$ be a CW-complex obtained by attaching three $4$-cells to $M_g$ along the three attaching regions of connected sum (see Figure \ref{CWcpx}). 

\begin{figure}[htbp]
\begin{center}
\includegraphics[width=130mm]{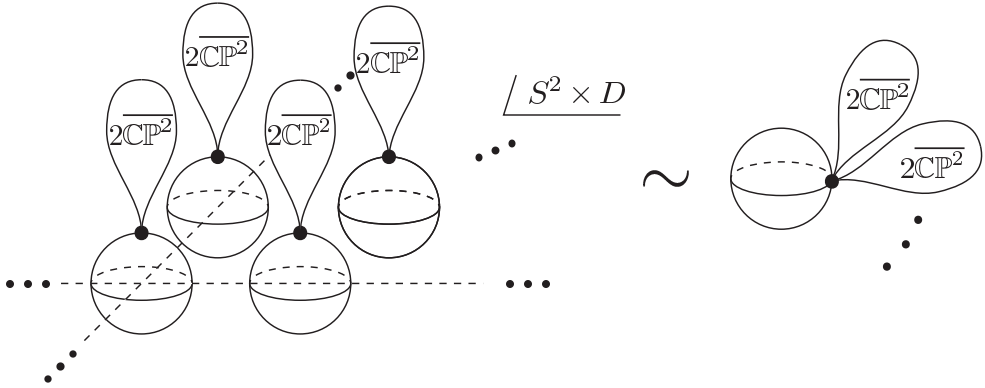}
\caption{Left: the figure describing $W_g$. 
Right: the wedge sum of $S^2$ and countably many $\overline{\mathbb{CP}^2}\vee\overline{\mathbb{CP}^2}$, 
which is obtained by collapsing $D$ to a point. }
\label{cover}
\end{center}
\end{figure}

Let $\iota:M_g\rightarrow X_g$ be the natural inclusion. 
By the cellular approximation theorem (for this theorem, see \cite{Ha}), the following map is isomorphism: 
\[
\iota_\ast:\pi_2(M_g,p_0)\rightarrow \pi_2(X_g,p_0). 
\]
Since $X_g$ is homotopic to $S^2\times \Sigma_{g-1}\vee S^1\times S^3\vee\overline{\mathbb{CP}^2}\vee\overline{\mathbb{CP}^2}$, 
we obtain: 
\begin{align*}
\pi_2(X_g,p_0) & \cong \pi_2(S^2\times \Sigma_{g-1}\vee S^1\times S^3\vee\overline{\mathbb{CP}^2}\vee\overline{\mathbb{CP}^2},p_0) \\
& \cong \pi_2(S^2\times \Sigma_{g-1}\vee S^1\vee S^3\vee\overline{\mathbb{CP}^2}\vee\overline{\mathbb{CP}^2},p_0), 
\end{align*}
where the second isomorphism is obtained by the cellular approximation theorem. 
We put $Y_g=S^2\times \Sigma_{g-1}\vee D^2 \vee S^3\vee\overline{\mathbb{CP}^2}\vee\overline{\mathbb{CP}^2}$ 
and denote by $j:X_g\rightarrow Y_g$ the inclusion map. 
Since $D^2$ is contractible and $S^3$ consists of the $0$-cell and the $3$-cell, we obtain: 
\[
\pi_2(Y_g,p_0)\cong \pi_2(S^2\times \Sigma_{g-1}\vee\overline{\mathbb{CP}^2}\vee\overline{\mathbb{CP}^2},p_0). 
\]
We denote by $W_g$ the universal cover of $S^2\times \Sigma_{g-1}\vee\overline{\mathbb{CP}^2}\vee\overline{\mathbb{CP}^2}$. 
$W_g$ is obtained by attaching countably many $\overline{\mathbb{CP}^2}\vee\overline{\mathbb{CP}^2}$ to $S^2\times D$, where $D$ is the universal cover of the closed surface, 
and is homotopic to $S^2\hspace{.85em}\bigvee\hspace{-2.3em}\raisebox{-1em}{\footnotesize{$\mu\in\pi_1(\Sigma_{g-1},q_0)$}}(\overline{\mathbb{CP}^2}\vee\overline{\mathbb{CP}^2})_{\mu}$ (see Figure \ref{cover}).

In general, the second homotopy group of a CW-complex is isomorphic to that of the universal cover of the complex. 
Thus, we obtain:
{\allowdisplaybreaks 
\begin{align*}
\pi_2(S^2\times \Sigma_{g-1}\vee\overline{\mathbb{CP}^2}\vee\overline{\mathbb{CP}^2},p_0) & \cong
\pi_2(S^2\hspace{.85em}\bigvee\hspace{-2.3em}\raisebox{-1em}{\footnotesize{$\mu\in\pi_1(\Sigma_{g-1},q_0)$}}(\overline{\mathbb{CP}^2}\vee\overline{\mathbb{CP}^2})_{\mu},p_0) \\
& \cong \mathbb{Z} \hspace{.85em}\bigoplus\hspace{-2.3em}\raisebox{-1em}{\footnotesize{$\mu\in\pi_1(\Sigma_{g-1},q_0)$}}(\mathbb{Z}\oplus\mathbb{Z})_{\mu}. 
\end{align*}
}

Eventually, we obtain the following group homomorphism: 
\[
\Phi=p_\ast^{-1}\circ j_\ast \circ \iota_\ast : \pi_2(M_g,p_0)\rightarrow 
\mathbb{Z} \hspace{.85em}\bigoplus\hspace{-2.3em}\raisebox{-1em}{\footnotesize{$\mu\in\pi_1(\Sigma_{g-1},q_0)$}}(\mathbb{Z}\oplus\mathbb{Z})_{\mu}. 
\]
If $\Phi(s)=(l,(m_{\mu},n_{\mu})_{\mu})\in\mathbb{Z} \hspace{.85em}\bigoplus\hspace{-2.3em}\raisebox{-1em}{\footnotesize{$\mu\in\pi_1(\Sigma_{g-1},q_0)$}}(\mathbb{Z}\oplus\mathbb{Z})_{\mu}$ for an element $s\in \pi_2(M_g,p_0)$, 
we have $\Phi(\gamma\cdot s)=(l,(m_{\mu},n_{\mu})_{\lambda\cdot\mu})$ for an element $\gamma=(\lambda,z)\in\pi_1(\Sigma_{g-1},q_0)\oplus\mathbb{Z}\cong\pi_1(M_g,p_0)$.

\begin{figure}[htbp]
\begin{center}
\includegraphics[width=135mm]{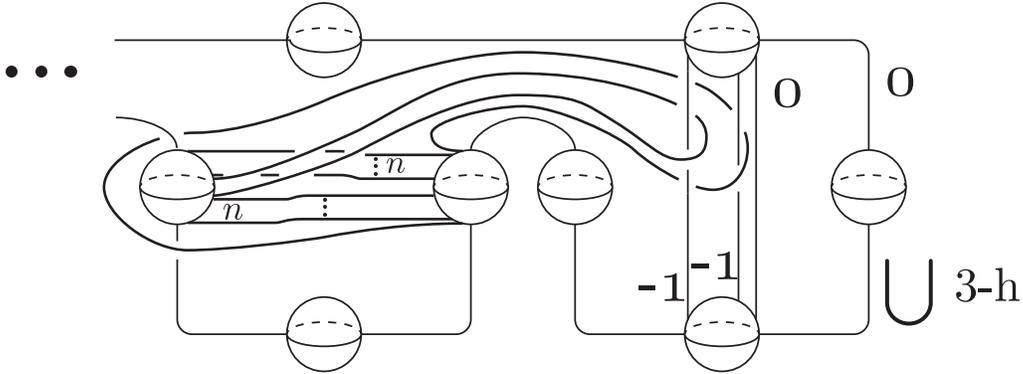}
\caption{the bold curve represents the section $\sigma_n$. }
\label{locus}
\end{center}
\end{figure}

For an integer $n$, let $\sigma_n:S^2\rightarrow M_g$ be a section whose image intersects the boundary of the lower side of $f_g$ at the locus illustrated in Figure \ref{locus}. 
Such a section exists since we can trivialize the locus illustrated in Figure \ref{locus} in the boundary of a regular neighborhood of a regular fiber in the higher and lower side of $f_g$. 

\begin{figure}[htbp]
\begin{center}
\includegraphics[width=45mm]{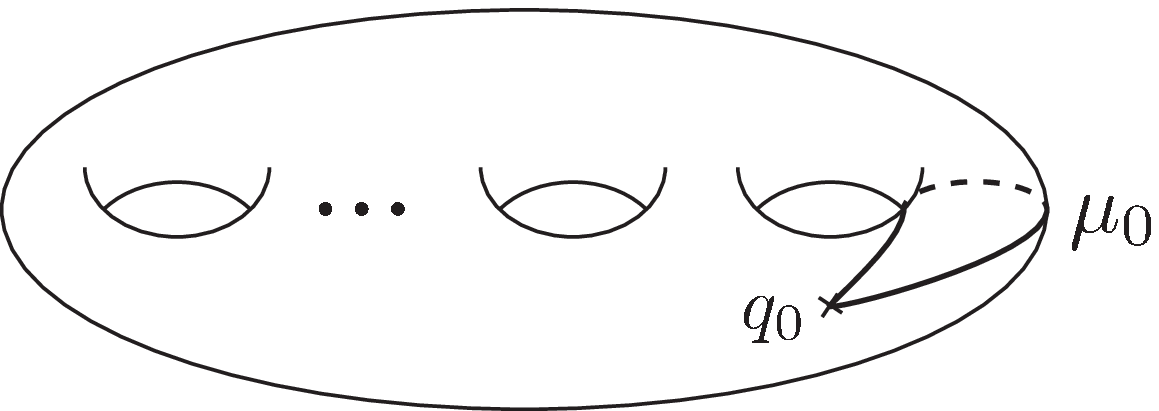}
\caption{}
\label{locuselement}
\end{center}
\end{figure}

We assume that the image of $\sigma_n$ contains the base point $p_0$, 
and regard $\sigma_n$ as an element in $\pi_2(M_g,p_0)$. 
By construction of $\sigma_n$, there exists an element $\gamma_n$ 
such that $\Phi(\gamma_n\cdot \sigma_n)$ is equal to $(1,(\delta_{1,\mu},\delta_{{\mu_0}^n,\mu})_{\mu})$, 
where $\delta_{\nu,\mu}$ is equal to $1$ if $\nu=\mu$ and $0$ otherwise, 
and $\mu_0\in\pi_1(\Sigma_{g-1},q_0)$ is the element described in Figure \ref{locuselement}


If $n$ is not equal to $m$, $(1,(\delta_{1,\mu},\delta_{{\mu_0}^n,\mu})_{\mu})$ is not equal to $(1,(\delta_{1,\mu},\delta_{{\mu_0}^m,\mu})_{\lambda\cdot\mu})$ for any elements $\lambda\in\pi_1(\Sigma_{g-1},q_0)\setminus \{1\}$. 
This means that $\sigma_n$ is not homotopic to $\sigma_m$ if $n$ is not equal to $m$. 
This completes the proof of Theorem \ref{main1}. \hfill $\square$

\section{Self-intersection of sections}\label{pro2}

We denote by $\Sigma_{g,1}$ the compact oriented surface with connected boundary 
and by $\delta$ a simple closed curve in $\Sigma_{g,1}$ parallel to the boundary. 
Let $\mathcal{M}_{g,1}$ be the mapping class group of $\Sigma_{g,1}$. 
It is known that there exists the natural surjective homomorphism $\psi:\mathcal{M}_{g,1}\rightarrow \mathcal{M}_g$ induced by the inclusion map $i:\Sigma_{g,1}\rightarrow \Sigma_g$. 
For a non-separating simple closed curve $\tilde{c}$ in $\Sigma_{g,1}$, 
we define $\mathcal{M}_{g,1}(\tilde{c})$ and $\tilde{\varphi}_{\tilde{c}}:\mathcal{M}_{g,1}(\tilde{c})\rightarrow \mathcal{M}_{g-1,1}$ 
as we define $\mathcal{M}_{g}(c)$ and $\varphi_c$.

\begin{lem}\label{section}

Let $\tilde{d},\tilde{d}_1,\ldots,\tilde{d}_n$ be simple closed curves in $\Sigma_{g,1}$. 
Suppose that these simple closed curves satisfy the following conditions: 

\begin{enumerate}

\item $\tilde{d}$ is non-separating; 

\item $t_{\tilde{d}_1}\cdot\cdots\cdot t_{\tilde{d}_n}\in\mathcal{M}_{g,1}(\tilde{d})$; 

\item $\varphi_{\tilde{d}}(t_{\tilde{d}_1}\cdot\cdots\cdot t_{\tilde{d}_n})={t_{\delta}}^k$, for some integer $k$. 

\end{enumerate}
Then there exists a simplified BLF $f:M\rightarrow S^2$ such that $f$ has a section $\sigma$ with $\sigma^2=-k$. 

\end{lem}

\noindent
{\it (Proof)}: We prove this lemma by constructing an explicit simplified BLF satisfying the desired condition. 
We take a $2$-disk $D$ in $\Sigma_g$ and we identify $\Sigma_{g,1}$ with $\Sigma_g\setminus \text{int}D$. 
We denote by $A$ the collar neighborhood of $\partial\Sigma_{g,1}$ in $\Sigma_{g,1}$. 
We fix an identification $D\cong D^2$ and $A\cong S^1\times [1,2]$ so that $\partial\Sigma_{g,1}$ corresponds to $S^1\times\{1\}$ in $A$. 
Then the map ${t_{\delta}}^k$ is represented by the following map: 
\begin{align*}
x &\mapsto\begin{cases}
x & (x\in\Sigma_{g,1}\setminus A), \\
(\text{exp}(\sqrt{-1}\theta+2\pi k\sqrt{-1}(2-s)),s) & (x=(\text{exp}(\sqrt{-1}\theta),s)\in A\cong S^1\times [1,2]). 
\end{cases}
\end{align*}

We first construct an LF over $D^2$ by attaching $n$ $2$-handles to $D^2\times\Sigma_g$
along $i(\tilde{d}_1),\ldots,i(\tilde{d}_n)$ in a regular fiber of $S^1\times\Sigma_g\subset D^2\times\Sigma_g$ 
with framing $-1$ with respect to the framing along a regular fiber (such a construction was introduced by Kas \cite{Kas}). 
Since $t_{\tilde{d}_1}\cdot\cdots\cdot t_{\tilde{d}_n}\in\mathcal{M}_{g,1}(\tilde{d})$, 
we can obtain a BLF over $D^2$ by round $2$-handle attachment (for details about this construction, see \cite{Ba}). 

By the condition (3) in the statement, the boundary of the resulting BLF is described as follows: 
\[
\Sigma_{g,1}\times I/((x,1)\sim({t_{\delta}}^k(x),0)) \cup D\times I/((x,1)\sim (x,0)). 
\]
Moreover, this BLF has a section $\tilde{\sigma}$ whose boundary is $\{0\}\times I/((x,1)\sim(x,0))$, 
where $0\in D$ is the center of the $2$-disk. 

To obtain a simplified BLF, we attach the trivial bundle $\Sigma_g\times D^2$ to the above BLF by the map 
\[
\Phi:\Sigma_g\times I/((x,1)\sim (x,0))\rightarrow\Sigma_{g,1}\times I/((x,1)\sim ({t_{\delta}}^k(x),0)) \cup D\times I/((x,1)\sim (x,0))
\]
defined as follows: 
\begin{align*}
\Phi(x,t) & = \begin{cases}
(x,t) & (x\in\Sigma_{g,1}\setminus A), \\
((\text{exp}(\sqrt{-1}\theta+2\pi k\sqrt{-1}t(s-2)),s),t) & (x=(\text{exp}(\sqrt{-1}\theta),s)\in A), \\
(r\text{exp}(\sqrt{-1}\theta-2\pi k\sqrt{-1}t),t) & (x=r\text{exp}(\sqrt{-1}\theta)\in D). 
\end{cases}
\end{align*}
The resulting simplified BLF has a section $\sigma=\{0\}\times D^2\cup_{\Phi}\tilde{\sigma}$. 
By the construction, the self-intersection of $\sigma$ is equal to $-k$. 
This completes the proof of Lemma \ref{section}. \hfill $\square$

\begin{figure}[htbp]
\begin{center}
\includegraphics[width=100mm]{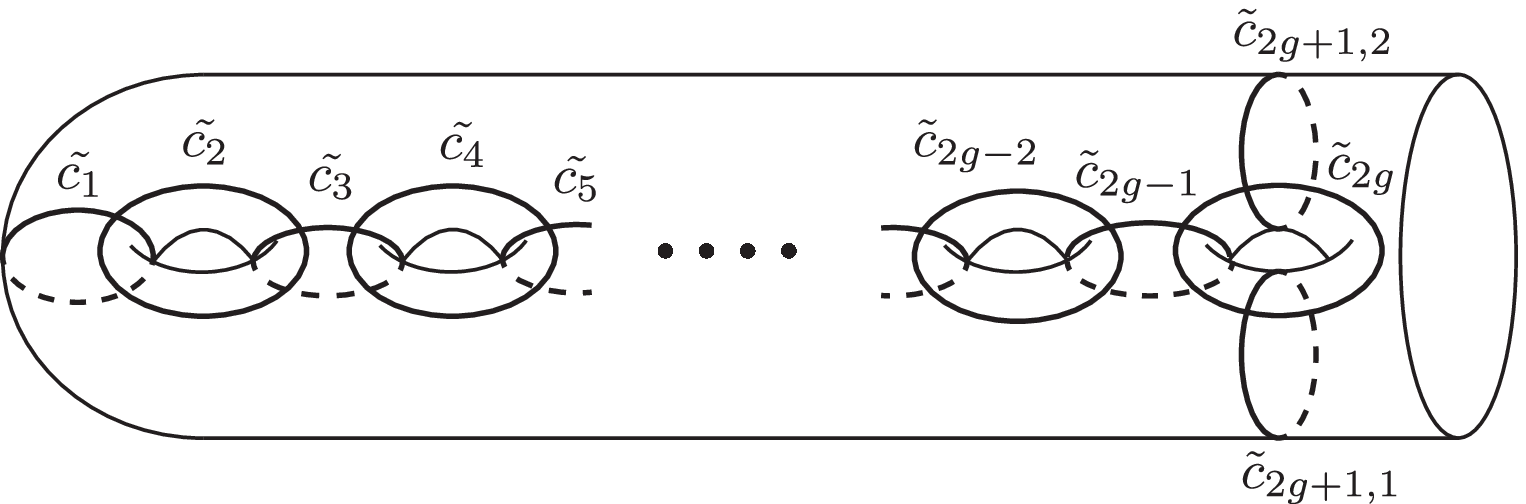}
\end{center}
\caption{}
\label{scc}
\end{figure}

\noindent
{\it (Proof of Theorem \ref{main2})}: 
We take simple closed curves $\tilde{c}_1,\ldots,\tilde{c}_{2g},\tilde{c}_{2g+1,1},\tilde{c}_{2g+1,2}$ in $\Sigma_{g,1}$ as shown in Figure \ref{scc}. 

\begin{figure}[htbp]
\begin{center}
\includegraphics[width=130mm]{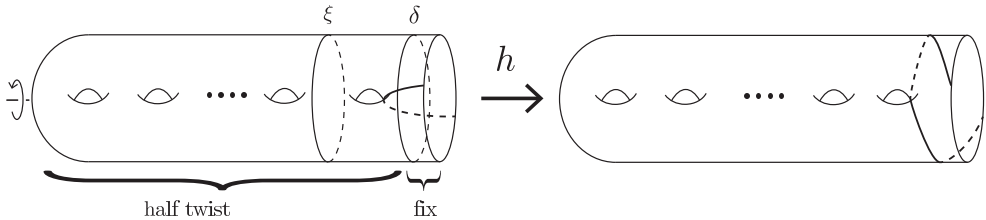}
\end{center}
\caption{$h$ twists the left side of the curve $\delta$ and fixes the right side of $\delta$ in the figure. }
\label{halftwist2}
\end{figure}

There exist the following relations in $\mathcal{M}_{g,1}$ for $g\geq 2$: 

\begin{enumerate}[(1)]

\item $(t_{\tilde{c}_1}\cdot\cdots\cdot t_{\tilde{c}_{2g-2}})^{2(2g-1)}=t_{\xi}$, 

\item $(t_{\tilde{c}_1}\cdot\cdots\cdot t_{\tilde{c}_{2g-1}})^{2g}=t_{\tilde{c}_{2g+1,1}}\cdot t_{\tilde{c}_{2g+1,2}}$, 

\item $(t_{\tilde{c}_1}\cdot\cdots\cdot t_{\tilde{c}_{2g}})^{2g+1}=h$, 

\end{enumerate}
where $\xi$ is the simple closed curve described in Figure \ref{halftwist2} and $h$ is the element of $\mathcal{M}_{g,1}$ as shown in Figure \ref{halftwist2}. 

By using these relations, we obtain the following relation: 
\[
t_{\tilde{c}_{2g}}\cdot\cdots\cdot t_{\tilde{c}_{2}}\cdot t_{\tilde{c}_{1}}^2\cdot t_{\tilde{c}_{2}}\cdot\cdots\cdot t_{\tilde{c}_{2g}}=t_{\tilde{c}_{2g+1,1}}^{-1}\cdot t_{\tilde{c}_{2g+1,2}}^{-1}\cdot h. 
\]
Since ${h}^2=t_{\delta}$, we obtain: 
\begin{align*}
(t_{\tilde{c}_1}\cdot\cdots\cdot t_{\tilde{c}_{2g}})^{(4g+2)n} & =t_{\delta}^n, \\
(t_{\tilde{c}_{2g}}\cdot\cdots\cdot t_{\tilde{c}_{2}}\cdot t_{\tilde{c}_{1}}^2\cdot t_{\tilde{c}_{2}}\cdot\cdots\cdot t_{\tilde{c}_{2g}})^{2n} & =t_{\tilde{c}_{2g+1,1}}^{-2n}\cdot t_{\tilde{c}_{2g+1,2}}^{-2n}\cdot t_{\delta}^n, \\
(t_{\tilde{c}_{2g}}\cdot\cdots\cdot t_{\tilde{c}_{2}}\cdot t_{\tilde{c}_{1}}^2\cdot t_{\tilde{c}_{2}}\cdot\cdots\cdot t_{\tilde{c}_{2g}})^{2n}\cdot(t_{\tilde{c}_1}\cdot\cdots\cdot t_{\tilde{c}_{2g-2}})^{2(2g-1)n} & =t_{\tilde{c}_{2g+1,1}}^{-2n}\cdot t_{\tilde{c}_{2g+1,2}}^{-2n}\cdot t_{\xi}^{n}\cdot t_{\delta}^{n}, 
\end{align*}
where $n$ is a positive integer. 
The right side of the above equations are in $\mathcal{M}_{g,1}(\tilde{c}_{2g+1,1})$. 
Since $\tilde{\varphi}_{\tilde{c}_{2g+1,1}}(t_{\tilde{c}_{2g+1,1}})=1$, 
$\tilde{\varphi}_{\tilde{c}_{2g+1,1}}(t_{\xi})=\tilde{\varphi}_{\tilde{c}_{2g+1,1}}(t_{\delta})=\tilde{\varphi}_{\tilde{c}_{2g+1,1}}(t_{\tilde{c}_{2g+1,2}})=t_{\delta}$ and 
$\tilde{\varphi}_{\tilde{c}_{2g+1,1}}(h)=h$, we obtain: 
\begin{align*}
\tilde{\varphi}_{\tilde{c}_{2g+1,1}}(t_{\delta}^n) & =t_{\delta}^n, \\
\tilde{\varphi}_{\tilde{c}_{2g+1,1}}(t_{\tilde{c}_{2g+1,1}}^{-2n}\cdot t_{\tilde{c}_{2g+1,2}}^{-2n}\cdot t_{\delta}^n) & =t_{\delta}^{-n}, \\
\tilde{\varphi}_{\tilde{c}_{2g+1,1}}(t_{\tilde{c}_{2g+1,1}}^{-2n}\cdot t_{\tilde{c}_{2g+1,2}}^{-2n}\cdot t_\xi^{n}\cdot t_\delta^{n}) & =1. 
\end{align*}
Thus, the conclusion holds by Lemma \ref{section}. \hfill $\square$

\section{Spin structures}\label{spin}

In this section, we discuss spin structures of total spaces of simplified BLFs. 

\par

Let $f:M\rightarrow S^2$ be a simplified BLF. 
Denote by $F\subset M$ a regular fiber in the lower side of $f$. 
A homology class $S\in H_2(M;\mathbb{Z})$ is called a {\it dual} of $F$ 
if the intersection number $S\cdot [F]$ is equal to $1$. 
If $f$ has a section $\sigma:S^2\rightarrow M$, the element $[\sigma(S^2)]$ is a dual of $F$. 
It is also easy to see that a dual of $F$ exists if the union of the higher side and the round cobordism of $f$ is simply connected. 

\begin{thm}\label{main3}

Let $f:M\rightarrow S^2$ be a genus-$g$ simplified BLF and $F$ a regular fiber in the lower side of $f$. 
We denote by $d_1,\ldots,d_n\subset\Sigma_g$ and $d\subset\Sigma_g$ vanishing cycles of Lefschetz singularities and the indefinite fold of $f$, respectively. 
Suppose that there exists a dual $S\in H_2(M;\mathbb{Z})$ of $F$. 
Then $M$ admits a spin structure if and only if the following two conditions hold: 

\begin{enumerate}[(a)]

\item there exists a quadratic form $q:H_1(\Sigma_g;\mathbb{Z}/2\mathbb{Z})\rightarrow\mathbb{Z}/2\mathbb{Z}$ with respect to the intersection form of $\Sigma_g$ 
such that $q(d)=0$ and $q(d_i)=1$ for all $i\in\{1,\ldots n\}$; 

\item the self-intersection of $S$ is even. 

\end{enumerate}

\end{thm}

\noindent
{\it (Proof)}: 
We first prove that the condition (a) in the statement holds if and only if 
the union of the higher side and the round cobordism of $f$ admits a spin structure. 
Let $\tilde{F}$ be a regular fiber in the higher side of $f$ and $\nu \tilde{F}\cong D^2\times\Sigma_g$ a regular neighborhood of $\tilde{F}$. 
It is known that $\nu \tilde{F}$ may admits exactly $2^{2g}$ distinct spin structures and that 
there exists one to one correspondence between the set of equivalence classes of spin structures of $\nu \tilde{F}$ and 
the set of quadratic forms $q:H_1(\Sigma_g;\mathbb{Z}/2\mathbb{Z})\rightarrow \mathbb{Z}/2\mathbb{Z}$. 
For a given spin structure $s$ of $\nu \tilde{F}$, the corresponding quadratic form $q_s$ is defined as follows: 
for an element $\gamma\in H_1(\Sigma_g;\mathbb{Z}/2\mathbb{Z})$, we take a simple closed curve $c\subset\Sigma_g\cong \tilde{F}$ 
which represents $\gamma$. 
Then $q_s(\gamma)$ is equal to $0$ if the restriction of $s$ to $c$ can be extended to the spin structure of the $2$-disk whose boundary is $c$ 
and is equal to $1$ otherwise (the reader should turn to \cite{Sti2} for more details about this correspondence). 

By the argument in \cite{Sti2}, the higher side of $f$ admits a spin structure if and only if 
there exists the quadratic form $q:H_1(\Sigma_g;\mathbb{Z}/2\mathbb{Z})\rightarrow \mathbb{Z}/2\mathbb{Z}$ such that
$q(d_i)=1$ for all $i=1,\ldots,n$. 

By Lemma \ref{roundhandle}, the union of the higher side and the round cobordism of $f$ is obtained by attaching a $2$-handle and a $3$-handle to the higher side. 
Moreover, it is easy to see that the attaching map of the $2$-handle preserves the spin structure obtained by restricting $s$ to the simple closed curve $d$. 
So we can extend $s$ to the $2$-handle if and only if $s|_d$ can be extended to the bounding $2$-disk. 
Since the attaching region of the $3$-handle is diffeomorphic to $S^2\times D^1$ and has the unique spin structure, 
we can extend $s$ to the round cobordism of $f$ if and only if $q(d)=0$. 
So the condition (a) is equivalent to the condition that 
the union of the higher side and the round cobordism of $f$ admits a spin structure. 

Now we are ready to prove Theorem \ref{main3}. 
If $M$ admits a spin structure, then the union of the higher side and the round cobordism of $f$ also admits a spin structure. 
So the condition (a) holds. 
Since the intersection form of $M$ is even, the self-intersection of $S$ must be even. 

The converse direction is easily proved by the same argument as in \cite{Sti2}. \hfill $\square$






\begin{rem}

In \cite{Wil2}, Williams introduced a {\it surface diagram} $(\Sigma_g,\Gamma)$ of a $4$-manifold $M$, 
where $g\geq 3$ and $\Gamma=(\gamma_1,\ldots,\gamma_k)$ is a $\mathbb{Z}/k\mathbb{Z}$-indexed collection of simple closed curves in $\Sigma_g$. 
This diagram is defined by using a {\it simplified purely wrinkled fibration} $f:M\rightarrow S^2$. 
The simple closed curves in $\Gamma$ represent the vanishing cycles of indefinite fold of $f$ (for more details, see \cite{Wil2}). 

By using the modification defined by Lekili \cite{Lek}, 
we can change indefinite cusps into Lefschetz singularities and indefinite folds and 
we obtain the simplified BLF $h:M\rightarrow S^2$ from a simplified purely wrinkled fibration $f$ with the surface diagram $(\Sigma_g,\Gamma=(\gamma_1,\ldots,\gamma_k))$. 
Let $W_h=(d_1,\ldots,d_k)$ be the Hurwitz system of $h$, then the class $[d_i]\in H_1(\Sigma_g;\mathbb{Z}/2\mathbb{Z})$ is equal to $[\gamma_i]+[\gamma_{i+1}]$. 
For any quadratic form $q:H_1(\Sigma_g;\mathbb{Z}/2\mathbb{Z})\rightarrow \mathbb{Z}/2\mathbb{Z}$, 
the following equation holds: 
{\allowdisplaybreaks
\begin{align*}
q([d_i]) &=q([\gamma_i])+q([\gamma_{i+1}])+[\gamma_i]\cdot[\gamma_{i+1}] \\
&=q([\gamma_i])+q([\gamma_{i+1}])+1. 
\end{align*}
}
So $q(d_i)$ is equal to $1$ if and only if $q(\gamma_i)=q(\gamma_{i+1})$. 
Thus, we obtain the following corollary. 

\begin{cor}

Let $f:M\rightarrow S^2$ be a simplified purely wrinkled fibration and 
$(\Sigma_g,\Gamma)$ a surface diagram of $M$ induced by $f$. 
Denote by $F$ a regular fiber of the lower side of $f$. 
We assume that there exists a dual $S\in H_2(M;\mathbb{Z})$ of $F$. 
Then $M$ admits a spin structure if and only if the following conditions hold: 

\begin{enumerate}[(a)]

\item there exists a quadratic form $q:H_1(\Sigma_g;\mathbb{Z}/2\mathbb{Z})\rightarrow \mathbb{Z}/2\mathbb{Z}$ 
such that $q(d)=0$ for all $d\in \Gamma$; 

\item the self-intersection of $S$ is even. 

\end{enumerate}

\end{cor}

\end{rem}

In the rest of this section, we will give some applications of Theorem \ref{main3}. 

\begin{exmp}

For an integer $n$ and a positive even integer $g=2k$, we denote by $f_{g,n}:M\rightarrow S^2$ the genus-$g$ simplified BLF constructed in the proof of Theorem \ref{main2} as a fibration with a section of square $n$. 
The Hurwitz system of $f_{g,n}$ is given as follows: 
\begin{align*}
(t_{{c}_1}\cdot\cdots\cdot t_{{c}_{2g}})^{(4g+2)|n|} & \hspace{1em}\text{(if $n$ is negative)}, \\
(t_{{c}_{2g}}\cdot\cdots\cdot t_{{c}_{2}}\cdot t_{{c}_{1}}^2\cdot t_{{c}_{2}}\cdot\cdots\cdot t_{{c}_{2g}})^{2}\cdot(t_{{c}_1}\cdot\cdots\cdot t_{{c}_{2g-2}})^{2(2g-1)} & \hspace{1em}\text{(if $n$ is zero)}, \\
(t_{{c}_{2g}}\cdot\cdots\cdot t_{{c}_{2}}\cdot t_{{c}_{1}}^2\cdot t_{{c}_{2}}\cdot\cdots\cdot t_{{c}_{2g}})^{2n} & \hspace{1em}\text{(if $n$ is positive)}, \\
\end{align*}

\noindent
where the simple closed curves $c_1,\ldots,c_{2g+1}$ is described in Figure \ref{scc2}. 
The group $H_1(\Sigma_g;\mathbb{Z}/2\mathbb{Z})$ is generated by the elements $[\gamma_1],\ldots,[\gamma_{2g}]$, 
where $\gamma_1,\ldots,\gamma_{2g}\subset\Sigma_g$ is simple closed curves described in Figure \ref{scc2}. 

\begin{figure}[htbp]
\begin{center}
\includegraphics[width=90mm]{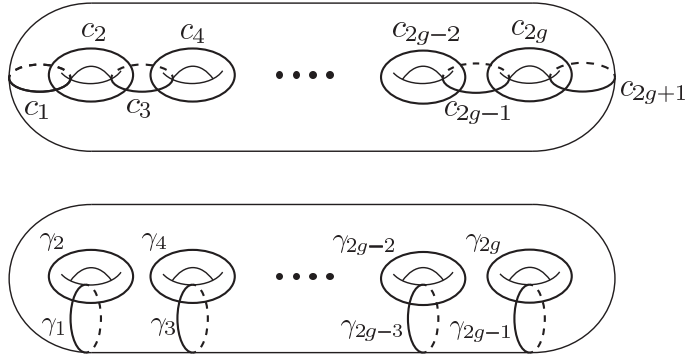}
\end{center}
\caption{simple closed curves on $\Sigma_g$}
\label{scc2}
\end{figure}

Let $q:H_1(\Sigma_g;\mathbb{Z}/2\mathbb{Z})\rightarrow \mathbb{Z}/2\mathbb{Z}$ be the quadratic form 
with respect to the intersection form of $\Sigma_g$ such that 
$q([\gamma_{2i}])=1$ for all $i=1,\ldots,g$, $q([\gamma_{4j-3}])=1$ and $q([\gamma_{4j-1}])=0$ for all $j=1,\ldots,k$. 
Since $[c_1]=[\gamma_1]$, $[c_{2g+1}]=[\gamma_{2g-1}]$, $[c_{2i}]=[\gamma_{2i}]$ ($i=1,\ldots,g$) and $[c_{2j+1}]=[\gamma_{2j-1}]+[\gamma_{2j+1}]$ ($j=1,\ldots,g-1$), 
we can calculate the value $q([c_i])$ as follows: 

{\allowdisplaybreaks
\begin{align*}
q([c_1]) & =q([\gamma_1])=1, \\
q([c_{2g+1}]) & =q([\gamma_{4k-1}])=0, \\
q([c_{2i}]) & =q([\gamma_{2i}])=1 \hspace{1em}\text{($i=1,\ldots,g$)}, \\
q([c_{2j+1}]) & =q([\gamma_{2j-1}])+q([\gamma_{2j+1}])+[\gamma_{2j-1}]\cdot[\gamma_{2j+1}] \\
& =1+0=1 \hspace{1em}\text{($j=1,\ldots,g-1$)}. 
\end{align*}
}

\noindent
So $q$ satisfies the condition (a) of (ii) in Theorem \ref{main3} for $f_{g,n}$. 
Moreover, $f_{g,n}$ has a section of square $n$. 
Thus, the total space of $f_{g,n}$ admits a spin structure if $n$ is even. 

\end{exmp}

We can completely classify {\it spin} genus-$1$ simplified BLF. 

\begin{prop}\label{clsspin}

Let $f:M\rightarrow S^2$ be a genus-$1$ simplified BLF. 
We assume that $M$ admits a spin structure and that $f$ has both Lefschetz and indefinite fold singularities. 
Then $M$ is diffeomorphic to $\sharp kS^2\times S^2$ for some $k\geq 1$. 

\end{prop}

\begin{rem}

Moishezon and Kas completely classified genus-$1$ LFs over $S^2$. 
Baykur, Kamada \cite{BK} and the author \cite{H} classified genus-$1$ simplified BLFs without Lefschetz singularities. 
So all we need to consider is the case $f$ has both Lefschetz and indefinite fold singularities. 

\end{rem}

\noindent
{\it (Proof of Proposition \ref{clsspin})}: 
We take simple closed curves $c_1,c_2\subset T^2$ so that the class $[c_1],[c_2]\in H_1(T^2;\mathbb{Z})$ is a generator of $H_1(T^2;\mathbb{Z})$ 
and that $c_1\cdot c_2=1$. 
Denote by $X_i\in \mathcal{M}_1=\text{Diff}^+(T^2)/\text{(isotopy)}$ ($i=1,2$) the right-handed Dehn twist along $c_i$. 
When we identify $\mathcal{M}_1$ with $SL(2,\mathbb{Z})$ by a suitable isomorphism, 
$X_1$ and $X_2$ correspond to the matrices 
$\begin{pmatrix}
1 & 0 \\
1 & 1 
\end{pmatrix}
$ and 
$\begin{pmatrix}
1 & -1 \\
0 & 1 
\end{pmatrix}
$, respectively. 
We define the sequences of elements of $SL(2,\mathbb{Z})$ $S_r$ and $T(n_1,\ldots,n_s)$ as follows: 
\begin{align*}
S_r &= (X_1,\ldots,X_1) \hspace{.5em} \text{($r$ $X_1$ stand in a line. )}, \\
T(n_1,\ldots,n_s) &=({X_1}^{-n_1}X_2{X_1}^{n_1},\ldots,{X_1}^{-n_s}X_2{X_1}^{n_s}). 
\end{align*}
By Theorem 3.11 in \cite{H}, we can assume that a Hurwitz system $W_f$ of $f$ is equal to $S_r\cdot T(n_1,\ldots,n_s)$ for some $r,n_1,\ldots,n_s\in\mathbb{Z}$
and that $w(W_f)$ corresponds to $\pm {X_1}^n$ for $n\in\mathbb{Z}$, 
where $w(W_f)\in SL(2,\mathbb{Z})$ is the product of all elements in $W_f$

If $r$ were not equal to $0$, $M$ would contain $\overline{\mathbb{CP}^2}$ as a connected sum component 
and $M$ would not be spin. 
So we have $r=0$. 

The vanishing cycles of Lefschetz singularities (resp. indefinite fold) of $f$ are $t_{c_1}^{n_1}(c_2),\ldots,t_{c_1}^{n_s}(c_2)$ (resp. $c_1$). 
By the Picard-Lefschetz formula, we obtain: 
\[
[t_{c_1}^{n_i}(c_2)]=[c_2]+n_i[c_1]\in H_1(T^2;\mathbb{Z}). 
\]

By Theorem \ref{main3}, there exists a quadratic form $q:H_1(T^2;\mathbb{Z}/2\mathbb{Z})\rightarrow \mathbb{Z}/2\mathbb{Z}$ 
such that $q([c_1])=0$ and $q([c_2]+n_i[c_1])=1$. 
On the other hand, there exists exactly two quadratic forms $q_0$, $q_1$ which satisfy $q_j([c_1])=0$ ($q_0([c_2])=0$, while $q_1([c_2])=1$). 
$q_j([c_2]+n_i[c_1])$ is calculated as follows: 
{\allowdisplaybreaks
\begin{align*}
q_j([c_2]+n_i[c_1]) &=\begin{cases}
q_j([c_2]) & \text{(if $n_i$ is even)}, \\
q_j([c_1]+[c_2])=q_j([c_2])+1 & \text{(if $n_i$ is odd)}, 
\end{cases} \\
& = \begin{cases}
0 & \text{(if $n_i$:even, $j=0$ or $n_i$:odd, $j=1$)}, \\
1 & \text{(if $n_i$:odd, $j=0$ of $n_i$:even, $j=1$)}. 
\end{cases}
\end{align*}
}

\noindent
Eventually, the integers $n_1,\ldots,n_s$ have same parity. 
In particular, the integers $n_1-n_2,\ldots,n_{s-1}-n_s$ are all even. 

It is known that the group $PSL(2,\mathbb{Z})$ has the following presentation: 
\[
PSL(2,\mathbb{Z})=<a,b|a^3,b^2>\cong\mathbb{Z}/3\mathbb{Z}\ast\mathbb{Z}/2\mathbb{Z}. 
\]
Let $p:SL(2,\mathbb{Z})\rightarrow PSL(2,\mathbb{Z})$ be the natural projection. 
Then $x_1=p(X_1)=aba$ and $x_2=p(X_2)=ba^2$. 
Since $w(W_f)=\pm {X_1}^n$, we obtain: 
{\allowdisplaybreaks
\begin{align*}
& {X_1}^{-n_1}X_2{X_1}^{n_1-n_2}\cdot\cdots\cdot{X_1}^{n_{s-1}-n_s}X_2{X_1}^{n_s}=\pm {X_1}^n, \\
\Rightarrow & {x_1}^{-n_1}x_2{x_1}^{n_1-n_2}\cdot\cdots\cdot{x_1}^{n_{s-1}-n_s}x_2{x_1}^{n_s}={x_1}^n, \\
\Rightarrow & x_2{x_1}^{n_1-n_2}\cdot\cdots\cdot{x_1}^{n_{s-1}-n_s}x_2={x_1}^m, 
\end{align*}
}

\noindent
where $m=n+n_1-n_s$. 

\begin{lem}\label{key}

Suppose that $n_i-n_{i+1}\neq 2$ for all $i\in\{1,\ldots,s-1\}$. 
Then $x_2{x_1}^{n_1-n_2}\cdot\cdots\cdot{x_1}^{n_{s-1}-n_s}x_2$ is equal to $bS$ or $a^2ba^2bS$, 
where $S=w_1\cdot\cdots\cdot w_k$ and $(w_1,\ldots,w_k)$ is a reduced sequence (i.e. $\{w_i,w_{i+1}\}=\{a,b\}\text{ or }\{a^2,b\}$) 
such that $w_1=a\text{ or }a^2$. 

\end{lem}

\noindent
{\it (Proof of Lemma \ref{key})}: 
We prove this statement by induction on $s$. 

We first look at the case $s=2$. 
$x_2{x_1}^{n_1-n_2}x_2$ is calculated as follows: 
{\allowdisplaybreaks
\begin{align*}
x_2{x_1}^{n_1-n_2}x_2 & = \begin{cases}
ba^2\cdot a(ba^2)^{n_1-n_2-1}ba\cdot ba^2 & \text{(if $n_1-n_2\geq 4$)}, \\
ba^2\cdot a^2(ba)^{-n_1+n_2-1}ba^2\cdot ba^2 & \text{(if $n_1-n_2\leq 0$)}, 
\end{cases} \\
& = \begin{cases}
a^2ba^2(ba^2)^{n_1-n_2-3}baba^2 & \text{(if $n_1-n_2\geq 4$)}, \\
(ba)^{-n_1+n_2}ba^2ba^2 & \text{(if $n_1-n_2\leq 0$)}. 
\end{cases}
\end{align*}
}
So the statement holds. 

We then look at the general case. 
By the induction hypothesis, we obtain: 
\[
x_2{x_1}^{n_2-n_3}\cdot\cdots\cdot{x_1}^{n_{s-1}-n_s}x_2=bS\text{ or }a^2ba^2bS, 
\]
where $S$ is the product of a reduced sequence starting from $a$ or $a^2$. 
We can calculate $x_2{x_1}^{n_1-n_2}$ as follows: 
{\allowdisplaybreaks
\begin{align*}
x_2{x_1}^{n_1-n_2} &= \begin{cases}
ba^2\cdot a(ba^2)^{n_1-n_2-1}ba & \text{(if $n_1-n_2\geq 4$)}, \\
ba^2\cdot a^2(ba)^{-n_1+n_2-1}ba^2 & \text{(if $n_1-n_2\leq 0$)}, 
\end{cases}\\
& = \begin{cases}
a^2ba^2(ba^2)^{n_1-n_2-4}ba^2ba & \text{(if $n_1-n_2\geq 4$)}, \\
(ba)^{-n_1+n_2}ba^2 & \text{(if $n_1-n_2\leq 0$)}, 
\end{cases}
\end{align*}
}
Hence, we obtain: 
{\allowdisplaybreaks
\begin{align*}
& x_2{x_1}^{n_1-n_2}\cdot\cdots\cdot{x_1}^{n_{s-1}-n_s}x_2 \\
= & \begin{cases}
a^2ba^2(ba^2)^{n_1-n_2-4}ba^2ba\cdot bS & \text{(if $n_1-n_2\geq 4$ and $x_2{x_1}^{n_2-n_3}\cdot\cdots\cdot{x_1}^{n_{s-1}-n_s}x_2=bS$)}, \\
a^2ba^2(ba^2)^{n_1-n_2-4}ba^2ba\cdot a^2ba^2bS & \text{(if $n_1-n_2\geq 4$ and $x_2{x_1}^{n_2-n_3}\cdot\cdots\cdot{x_1}^{n_{s-1}-n_s}x_2=a^2ba^2bS$)}, \\
(ba)^{-n_1+n_2}ba^2\cdot bS & \text{(if $n_1-n_2\leq 0$ and $x_2{x_1}^{n_2-n_3}\cdot\cdots\cdot{x_1}^{n_{s-1}-n_s}x_2=bS$)}, \\
(ba)^{-n_1+n_2}ba^2\cdot a^2ba^2bS  & \text{(if $n_1-n_2\leq 0$ and $x_2{x_1}^{n_2-n_3}\cdot\cdots\cdot{x_1}^{n_{s-1}-n_s}x_2=a^2ba^2bS$)}, 
\end{cases} \\
= & \begin{cases}
a^2ba^2(ba^2)^{n_1-n_2-4}ba^2babS & \text{(if $n_1-n_2\geq 4$ and $x_2{x_1}^{n_2-n_3}\cdot\cdots\cdot{x_1}^{n_{s-1}-n_s}x_2=bS$)}, \\
a^2ba^2(ba^2)^{n_1-n_2-4}babS & \text{(if $n_1-n_2\geq 4$ and $x_2{x_1}^{n_2-n_3}\cdot\cdots\cdot{x_1}^{n_{s-1}-n_s}x_2=a^2ba^2bS$)}, \\
(ba)^{-n_1+n_2}ba^2bS & \text{(if $n_1-n_2\leq 0$ and $x_2{x_1}^{n_2-n_3}\cdot\cdots\cdot{x_1}^{n_{s-1}-n_s}x_2=bS$)}, \\
(ba)^{-n_1+n_2}baba^2bS  & \text{(if $n_1-n_2\leq 0$ and $x_2{x_1}^{n_2-n_3}\cdot\cdots\cdot{x_1}^{n_{s-1}-n_s}x_2=a^2ba^2bS$)}. 
\end{cases}
\end{align*}
}
This completes the proof of Lemma \ref{key}. \hfill $\square$

By Lemma \ref{key}, $x_2x_1^{n_1-n_2}\cdots x_1^{n_{s-1}-n_s}x_2$ would not be equal to $x_1^m$ if $n_i-n_{i+1}\neq 2$ for all $i\in\{1,\ldots,s-1\}$.
Hence, we have $n_i-n_{i+1}=2$ for some $i\in\{1,\ldots,s-1\}$ and $M$ contains $S^2\times S^2$ as a connected sum component. 
By applying this argument successively, we can complete the proof of Proposition \ref{clsspin}. \hfill $\square$

\vspace{1em}

\noindent
{\bf Acknowledgments. }
The author would like to thank Hisaaki Endo for his helpful comments for the draft of this paper. 
The author also wishes to express his gratitude to Naoyuki Monden for his many useful suggestions, 
especially on self-intersection of sections. 
The author is supported by Yoshida Scholarship 'Master21' and he would like to thank Yoshida Scholarship Foundation for their support.

\end{document}